\documentclass{article}

\title{Regularity of $R(X)$ does not pass to finite unions}
\author{J. F. Feinstein}
\date{Creative Commons License: CC-BY-NC-ND}

 \usepackage{amsmath,amssymb,amsthm,amsfonts}
 \newtheorem{theorem}{Theorem}
  \newtheorem{corollary}[theorem]{Corollary}
    \newtheorem{lemma}[theorem]{Lemma}
    \newtheorem{proposition}[theorem]{Proposition}
\def\C{\mathbb{C}}
\def\N{\mathbb{N}}
\def\R{\mathbb{R}}
\def\Re{{\rm Re}}

\newcommand{\abs}[1]{\lvert #1 \rvert}

\newcommand{\Proof}{\noindent {\bf Proof. }}

\newcommand{\QED}{\hfill$\Box$}
\def\dist{\textrm{dist}}
\def\sskip{\vskip 0.3cm}

\begin{document}
\maketitle

\begin{abstract}
We show that there are compact plane sets $X$, $Y$ such that $R(X)$ and $R(Y)$
are regular but $R(X \cup Y)$ is not regular.
\end{abstract}

\section{Introduction}
In \cite{1963M}, McKissick gave the first example of a non-trivial, regular uniform algebra.
His example was $R(X)$ for a suitable compact plane set $X$.
(We recall the relevant definitions in the next section.)
In \cite[Corollary 7.7]{2016FMY} it was shown that if $X$ and $Y$ are compact plane sets such that
$R(X)$ and $R(Y)$ are regular and $X \cap Y$ is countable, then $R(X \cup Y)$ is also regular.

This note addresses the question of whether $R(X \cup Y)$ is regular whenever $R(X)$ and $R(Y)$
are regular, without any additional assumptions on the compact plane sets $X$ and $Y$.
If the answer was positive, the result would clearly also hold for all finite unions.
Bearing this in mind, we show that the answer is negative, in Theorem 6 below, by the slightly indirect method of finding
\textbf{four} compact plane sets $X_k$, $1 \leq k \leq 4$, such that each $R(X_k)$ is regular,
but such that $\displaystyle R\left(\bigcup_{k=1}^4 X_k\right)$ is not regular.

We describe this result by saying that \emph{regularity of $R(X)$ does not pass to finite unions}. Along the way, in Theorem 5, we give a slightly easier example showing that (in an obvious sense) regularity of $R(X)$ does not pass to countable unions.

\section{Preliminaries}

Throughout this note, by a \textit{compact space} we shall mean a non-empty, compact,
Hausdorff topological space;
by a \textit{compact plane set} we shall mean a {non-empty}, compact subset of the complex plane.
We shall use the term \emph{clopen} to describe sets which are both open and closed in a given topological space
(typically a compact plane set).

Let $a\in\C$ and let $r>0$.
We denote the open disk of radius $r$ and centre $a$ by $D(a,r)$ and the corresponding closed disk by $\bar D(a,r)$.
We denote the diameter of a non-empty, bounded subset $E$ of $\C$ by
diam$(E)$.

\sskip

We assume that the reader has some familiarity with uniform algebras.
We refer the reader to \cite{browder,Ga,St-book} for further background.
For the general theory of commutative Banach algebras, the reader may consult \cite{Dales,Kaniuth}.

\sskip

Let $X$ be a {compact space}, and let $C(X)$ be the algebra of all continuous complex-valued functions on
$X$. For each function $f \in C(X)$ and each non-empty subset $E$ of $X$, we denote the uniform norm of the
restriction of $f$ to $E$ by ${\abs {f}}_E$.
In particular, we denote by ${\abs {\,\cdot\,}}_X$ the uniform norm on $X$.
When endowed with the norm ${\abs {\,\cdot\,}}_X$, $C(X)$ is a Banach algebra.
A \emph{uniform algebra} on $X$ is a closed subalgebra of $C(X)$ that contains the
constant functions and separates
the points of $X$.
We say that a uniform algebra $A$ on $X$ is {\it nontrivial\/}
if $A\neq C(X)$, and is {\it natural\/} (on $X$)
if $X$ is the character space of $A$
(under the usual identification of points of $X$ with evaluation functionals).

Let $A$ be a natural uniform algebra on $X$, and let $x\in X$.
We denote by $J_x$ the ideal of functions $f$ in $A$ such that $x$ is
in the interior of the zero set of $f$, $f^{-1}(\{0\})$.
We denote by $M_x$ the ideal of functions $f$ in $A$ such that $f(x) = 0$.

We say that $x$ is a \emph{point of continuity} (for $A$) if, for all $y \in X \setminus \{x\}$
we have $ J_y\nsubseteq M_x$; %, there exists $f\in J_y$ with $f(x)\neq 0$
we say that $x$ is an \emph{R-point} (for $A$) if, for all $y \in X \setminus \{x\}$ we have
$J_x\nsubseteq M_y$. %, there exists $f\in J_x$ with  $f(y)\neq 0$.
We say that $A$ is \emph{regular} if, for every closed subset $F$ of $X$ and every
$y\in X\setminus F$, there exists $f\in A$ with $f(y)=1$ and $f(F)\subseteq \{0\}$.

The natural uniform algebra $A$ is regular
if and only if every point of $X$ is a point of continuity, and
this is also equivalent to the condition that every point of $X$ is an R-point
(\cite{2012FeinsteinMortini,FeinsteinSomerset2000}).

\sskip

Let $X$ be a compact plane set. By
$R(X)$ we denote the set of those functions $f\in C(X)$ which can be uniformly approximated on $X$
by rational functions with no poles on $X$.
It is standard that $R(X)$ is a natural uniform algebra on $X$.

Let $x \in X$. As we will sometimes need to work with different compact plane sets simultaneously, in this note
we will denote the ideals $J_x$ and $M_x$ in $R(X)$ by
$J^X_x$ and $M^X_x$ respectively.

We will occasionally use the fact that all idempotents in $C(X)$ are automatically in $R(X)$.
This does not require the full force of the Shilov Idempotent Theorem,
as for $R(X)$ it follows from Runge's theorem.

\sskip

Let $X$ be a compact plane set, let $Y$ be a compact subset of $X$ and suppose that $x,y \in Y$ with $x \neq y$.
Trivially the restriction $R(X)|_Y$ is contained in $R(Y)$. It follows immediately that if $J^X_x \nsubseteq M^X_y$ then $J^Y_x \nsubseteq M^Y_y$.
Thus if $x$ is an R-point for $R(X)$ then $x$ is also an R-point for $R(Y)$.
Similarly, if $x$ is a point of continuity for $R(X)$ then $x$ is a point of continuity for $R(Y)$.
If $R(X)$ is regular, then so is $R(Y)$.

\sskip

We now discuss some conditions under which converses to some of these implications hold.

\begin{lemma}
Let $X$ be a compact plane set and suppose that $E$ is a non-empty, clopen subset of $X$. Let $x \in E$.
If $x$ is an R-point for $R(E)$ then $x$ is an R-point for $R(X)$. Similarly, if $x$ is a point of continuity for $R(E)$ then $x$ is a point of continuity for $R(X)$.
\end{lemma}

\Proof
This is almost trivial in view of the idempotents available in $R(X)$. The only point that may be worth noting is the fact that, for all
$y \in E$, each function in $J_y^E$ has an extension in $J_y^X$. This is because every function in $R(E)$ can be extended to give a function in $R(X)$ which is constantly zero on $X\setminus E$.
\QED

\medskip

Let $X$ and $Y$ be compact plane sets, with $Y \subseteq X$, and suppose that $X$ contains no bounded component of $\C \setminus Y$.  Then, by Runge's theorem, $R(X)|_Y$ is dense in $R(Y)$.
This condition concerning the bounded components of $\C \setminus Y$ is, of course, automatically satisfied if $X$ has empty interior.

\begin{theorem}
Let $X$ and $Y$ be compact plane sets, with $Y \subseteq X$ and suppose that $x,y \in Y$ with $x \neq y$.
Suppose that $X$ contains no bounded component of $\C \setminus Y$ and that $X \setminus Y$ is the union of a sequence of pairwise disjoint, non-empty relatively clopen subsets $Y_n$ of $X$ whose diameters tend to $0$ as $n \to \infty$.
Then $Y$ is a peak set for $R(X)$ and $R(X)|_Y = R(Y)$. Moreover the following implications hold.

\begin{enumerate}
\item[(i)]If $J^Y_x  \nsubseteq M^Y_y$ then $J^X_x \nsubseteq M^X_y$.
\item[(ii)]If $x$ is an R-point for $R(Y)$ then $x$ is an R-point for $R(X)$.
\item[(iii)]If $x$ is a point of continuity for $R(Y)$ then $x$ is a point of continuity for $R(X)$.
\end{enumerate}
\end{theorem}

\Proof
In this proof, we work on $X$ when discussing clopen sets or characteristic functions: the characteristic functions of the clopen sets are then precisely the idempotents in $R(X)$.

Note that clopen subsets of $X$ are automatically peak sets (using idempotents in $R(X)$).
Clearly $Y$ is an intersection of (a sequence of) clopen subsets of $X$, so $Y$ is also a peak set for $R(X)$.
Thus $R(X)|_Y$ is a uniformly closed subset of $R(Y)$. As $X$ contains no bounded component of $\C \setminus Y$,
$R(X)|_Y$ is dense in $R(Y)$, so equality holds.

Using idempotents, we also see that, for each $z \in X \setminus Y$, we have $J^X_x \nsubseteq M^X_z$ and $J^X_z\nsubseteq M^X_x$. Thus (ii) and (iii) follow quickly from (i).

To prove (i), suppose that $J^Y_x  \nsubseteq M^Y_y$. Choose $f \in J^Y_x \setminus M^Y_y$. We shall show that $f$ has an extension in $J^X_x$, from which (i) follows. We know that $f$ has an extension, say $F$, in $R(X)$. However $F$ might not be in $J^X_x$, so we modify $F$ as follows.

Choose $r>0$ such that $f(\bar D(x,r) \cap Y) = \{0\}$.

Set
\[
S=\{n \in \N: Y_n \cap \bar D(x,r) \neq \emptyset\}\,.
\]
If $S$ is finite, then $F$ is already in $J^X_x$. So we may assume that $S$ is infinite, say
$S=\{n_1,n_2,\dots\}$ with $n_1<n_2<\cdots$.
Set $U = \bigcup_{n \in S} Y_n = \bigcup_{k=1}^\infty Y_{n_k}$.

We define $\tilde f : X \to \C$ by $\tilde f = (1-\chi_U) F$. We claim that $\tilde f$ is the extension of $f$ that we need. Clearly $\tilde f$ vanishes on $X \cap D(x,r)$. It remains to show that $\tilde f \in R(X)$ (and, in particular, that $\tilde f$ is continuous).

For each $m \in \N$, set $U_m=\bigcup_{k=1}^m Y_{n_k}$ and $V_m = \bigcup_{k=m+1}^\infty Y_{n_k}$. Then $U_m$ is clopen in $X$, so $\chi_{U_m} \in R(X)$ and hence $(1 - \chi_{U_m}) F  \in R(X)$. Since diam$(Y_{n_k}) \to 0$ as $k \to \infty$, we see that the the sequence of clopen sets $Y_{n_k}$ can accumulate only at points of $\bar D(x,r) \cap Y$, and that dist$(Y_{n_k},\bar D(x,r) \cap Y) \to 0$ as $k \to \infty$. It then follows that $|F|_{V_m} \to 0$ as $m \to \infty$, and hence
$(1 - \chi_{U_m}) F  \to \tilde f$ uniformly on $X$. The result follows.
\QED

\sskip

Using Theorem 2 and our earlier observations, we immediately obtain the following corollary concerning the regularity of $R(X)$ for this type of compact plane set $X$.

\begin{corollary}
Let $X$ and $Y$ be compact plane sets, with $Y \subseteq X$.
Suppose that $X$ contains no bounded component of $\C \setminus Y$ and that
$X \setminus Y$ is the union of a sequence of pairwise disjoint, non-empty relatively clopen subsets $Y_n$ of $X$ whose diameters tend to $0$ as $n \to \infty$.
Then $R(X)$ is regular if and only if $R(Y)$ is regular and each $R(Y_n)$ is regular.
\end{corollary}
\section{Main results}
We are almost ready to construct the examples mentioned in the introduction. First we need the following result, which is taken from \cite{2001Feinstein} (see also \cite{FM,FY}).

\begin{proposition}
There exist a compact plane set $X$ and a non-degenerate closed interval $J$ in $\R$ with $J\subseteq X$ such that $R(X)$ is not regular, but every point of $X\setminus J$ is a point of continuity for $R(X)$.
\end{proposition}

Let $X$, $J$ be as in Proposition 4. Clearly $X$ has empty interior.
For every non-empty compact subset $K$ of $X\setminus J$, $R(K)$ is regular, because every point of $K$ is a point of continuity for $R(K)$.
We have $R(J)=C(J)$, and so $R(J)$ is regular. However, $R(X)$ is not regular, so not every point of $X$ is a point of continuity for $R(X)$. The points of $X$ which are not points of continuity for $R(X)$ must all be in $J$. By Lemma 1, $J$ is not clopen in $X$. In particular, $X \neq J$.

\sskip

Using this example, we next show that, in a suitable sense, regularity of $R(X)$ does not pass to countable unions.

\begin{theorem}
There exist a compact plane set $X$ and a sequence of compact plane sets $X_k~(k \in \N)$ such that $X=\bigcup_{k\in \N} X_k$ and each $R(X_k)$ is regular, but $R(X)$ is not regular.
\end{theorem}

\textbf{Proof.}
Let $X$ and $J$ be as in Proposition 4. Scaling $X$ if necessary, we may assume that
$\{z\in X:\dist(z,J)\geq 1/2\} \neq \emptyset$.
Set $X_1=J$ and, for $k \geq 2$, set
$X_k=\{z\in X:\dist(z,J) \geq 1/k\,\}$.
Then $R(X)$ is not regular, but each $R(X_k)$ is regular, and $X=\bigcup_{k\in \N} X_k$ as required.
\QED
\sskip
Thus regularity of $R(X)$ does not pass to countable unions. We finish by showing that regularity of $R(X)$ does not pass to finite unions either.
\sskip

\begin{theorem}
There exist a compact plane set $X$ and four compact plane sets $X_k$, $k=1,2,3,4$, such that
$\displaystyle X=\bigcup_{k=1}^{4} X_k$ and each $R(X_k)$ is regular, but $R(X)$ is not regular.
\end{theorem}

\textbf{Proof.}
Again, let $X$ and $J$ be as in Proposition 4.
Scaling $X$ if necessary, we may assume that
$X \subseteq \{z\in \C:\dist(z,J)\leq 1\}$.

We begin by defining subsets $K_o$ and $K_e$ as follows.

Set
\[
K_o=J \cup \bigcup_{n\in\N} \{z \in X: 1/(2n+1) \leq \dist(z,J) \leq 1/2n\}
\]
and
\[
K_e=J \cup \bigcup_{n\in\N} \{z \in X: 1/2n \leq \dist(z,J) \leq 1/(2n-1)\}\,.
\]
Then $X = K_o \cup K_e$, and $R(X)$ is not regular.

It is not clear whether or not $R(K_o)$ or $R(K_e)$ are regular. In the following, we assume that $J$ is clopen in neither $K_o$ nor $K_e$ (otherwise the arguments may be simplified).  We now show that each of $K_o$ and $K_e$ is a union of two sets where we can establish regularity using Corollary 3. We may write $K_o=X_1 \cup X_3$ and $K_e=X_2 \cup X_4$, where the sets $X_k$ are defined with the help of vertical bands as follows:

\[
X_1=J \cup \bigcup_{n\in\N} \{z \in X: 1/(2n+1) \leq \dist(z,J) \leq 1/2n~\text{and}~\sin(n\Re(z))\leq 0 \}\,;
\]

\[
X_2=J \cup \bigcup_{n\in\N} \{z \in X: 1/2n \leq \dist(z,J) \leq 1/(2n-1)~\text{and}~\sin(n\Re(z))\leq 0 \}\,;
\]

\[
X_3=J \cup \bigcup_{n\in\N} \{z \in X: 1/(2n+1) \leq \dist(z,J) \leq 1/2n~\text{and}~\sin(n\Re(z))\geq 0 \}\,;
\]

\[
X_4=J \cup \bigcup_{n\in\N} \{z \in X: 1/2n \leq \dist(z,J) \leq 1/(2n-1)~\text{and}~\sin(n\Re(z))\geq 0 \}\,.
\]
 Then, for $1 \leq k \leq 4$, it is easy to see that $X_k$ is the disjoint union of $J$ with a sequence $(Y_{k,n})_{n=1}^\infty$ of pairwise disjoint, non-empty relatively clopen subsets of $X_k$ whose diameters tend to $0$. Each of the sets $X_k$ then satisfies the conditions of Corollary 3 with $Y=J$ and $Y_n=Y_{k,n}$ $(n\in \N)$. Since $R(J)$ and all of the $R(Y_{k,n})$ are regular, each $R(X_k)$ is regular by Corollary 3. Certainly $\displaystyle X=\bigcup_{k=1}^{4} X_k$, and $R(X)$ is not regular, as required.
\QED

\sskip
As noted earlier, it follows that regularity of $R(X)$ does not pass to unions of two sets either. However, we do not know which pair of sets mentioned above demonstrates this. Clearly at least one of the following pairs of sets must provide a suitable example: $X_1$ and $X_3$; $X_2$ and $X_4$; $K_o$ and $K_e$.

\sskip
We do not know whether or not regularity of $R(X)$ \emph{does} pass to unions of two sets if the sets are assumed to be connected.

\sskip
We are grateful to the referee for a careful reading and for useful comments and suggestions.

\bibliographystyle{amsplain}

\end{document}